\documentclass[12pt]{amsart}
\usepackage{amscd}

\newcommand{\h}{Hoch\-schild\ }
\newcommand{\hc}{Hoch\-schild co\-chain}
\newcommand{\hg}{ho\-mot\-o\-py G-}
\newcommand{\hl}{ho\-mot\-o\-py Lie\ }

\newcommand{\A}{\mathcal{A}}
\newcommand{\ass}{\mathcal{A}ssoc}
\newcommand{\BV}{\operatorname{BV}}
\newcommand{\coder}{\operatorname{Coder}}
\newcommand{\comm}{\mathcal{C}omm}

\newcommand{\End}[1]{\mathcal{E}nd\,_{#1}}

\newcommand{\F}[1]{\underline{\mathcal M}(#1)}
\newcommand{\FF}{\underline{\mathcal M}}
\newcommand{\FFR}{\underline{\mathcal M}_r}
\newcommand{\FR}[1]{\underline{\mathcal M}_r(#1)}
\newcommand{\G}{\mathcal{G}}
\newcommand{\gtg}{\mathfrak{g}}
\newcommand{\Hom}{\operatorname{Hom}}
\newcommand{\IN}{{\operatorname{in}}}
\newcommand{\lie}{\mathcal{L}ie}
\newcommand{\Mc}{\overline{\MM}}
\newcommand{\Mh}{\widehat{\mathcal{M}}}
\newcommand{\MM}{\mathcal{M}}
\newcommand{\Nc}{\overline{\mathcal{N}}}
\newcommand{\PP}{\mathcal{P}}
\newcommand{\QO}{\operatorname{QO}}
\newcommand{\Res}{\operatorname{Res}}
\newcommand{\Sym}{\mathcal{S}}
\newcommand{\Vr}{V_{\rm r}}
\newcommand{\X}[1]{\underline{\mathcal{M}}_{0,#1}}

\newcommand{\nc}{{\mathbb{C}}}
\newcommand{\nq}{{\mathbb{Q}}}
\newcommand{\nr}{{\mathbb{R}}}
\newcommand{\nz}{{\mathbb{Z}}}

\newtheorem{thm}{Theorem}[section]
\newtheorem{lm}[thm]{Lemma}
\newtheorem{prop}[thm]{Proposition}
\newtheorem{crl}[thm]{Corollary}
\newtheorem{conj}[thm]{Conjecture}
\newtheorem{quest}[thm]{Question}

\theoremstyle{definition}

\newtheorem{df}{Definition}[section]
\newtheorem{ex}{Example}[section]

\theoremstyle{remark}

\newtheorem{rem}{Remark}
\newtheorem{ack}{Acknowledgment}

\begin{document}

\title[Homotopy G-algebras and topological field theory]
{Homotopy Gerstenhaber algebras and topological field theory}

\author
[T. Kimura]{Takashi Kimura}
\address
{Department of Mathematics, Boston University, Boston, MA 02215}
\email{kimura@math.bu.edu}
\thanks{Research of the first author was supported in part by an NSF
postdoctoral research fellowship}

\author
[A. A. Voronov]{Alexander A. Voronov}
\address
{Department of Mathematics, University of Pennsylvania, Philadelphia, PA
19104-6395}
\email{voronov@math.upenn.edu}
\thanks{Research of the second author was supported in part by NSF grant
DMS-9402076}

\author
[G. J. Zuckerman]{Gregg J. Zuckerman}
\address
{Department of Mathematics, Yale University, New Haven, CT 06520}
\email{zuckerman@math.yale.edu}
\thanks{Research of the third author was supported in part by NSF grant
DMS-9307086}

\date{February 24, 1996}

\begin{abstract} We prove that the BRST complex of a topological conformal
field theory is a homotopy Gerstenhaber algebra, as conjectured by Lian
and Zuckerman in 1992. We also suggest a refinement of the original
conjecture for topological vertex operator algebras. We illustrate the
usefulness of our main tools, operads and ``string vertices'' by obtaining
new results on Vassiliev invariants of knots and double loop spaces.
\end{abstract}

\maketitle

        Two-dimensional topological quantum field theory (TQFT)
at its most elementary level is the theory of $\nz$-graded
commutative associative algebras (with some additional
structure) \cite{se}.  Thus, it came as something of a
surprise when several groups of mathematicians realized that
the physical state space of a 2D TQFT has the structure of a
$\nz$-graded Lie algebra, relative to a new grading equal to
the old grading minus one.  Moreover, the commutative and Lie
products fit together nicely to give the structure of a
Gerstenhaber algebra (G-algebra), a $\nz$-graded Poisson
algebra for which the Poisson bracket has degree $-1$ (see
Section \ref{g}).  This G-algebra structure is best understood
in the framework of 2D topological conformal field theories
(TCFTs) (see Section \ref{tcft}) wherein operads of moduli
spaces of Riemann surfaces play a fundamental role.

        G-algebras arose explicitly in M. Gerstenhaber's work
on the Hochschild cohomology theory for associative algebras
(see Section~\ref{g} for this and several other contexts for
the theory of G-algebras).  Operads arose in the work of J.
Stasheff, Gerstenhaber and later work of P. May on the
recognition problem for iterated loop spaces. Eventually, F.
Cohen discovered that the homology of a double loop space is
naturally a G-algebra, see Section~\ref{g}; in fact, a double
loop space is naturally an algebra over the little disks operad
employed by Cohen and also Boardman and Vogt, see Section~\ref{double}.
(The reader should consult the article
\cite{may} by May in these proceedings.)

        In joint work, B. Lian and G. Zuckerman (see also the
joint work of M. Penkava and A. Schwarz) discovered the
above-mentioned G-algebra structure in the context of
topological vertex operator algebras (TVOAs) (see Section
\ref{tvoa}), which are a powerful algebraic starting point for
the construction of 2D topological conformal field theories.
Lian and Zuckerman also gave a number of concrete constructions
of various examples of TVOAs, TCFTs and G-algebras.

        In an attempt to understand Lian-Zuckerman's work
geometrically, E. Getzler \cite{g} found a G-algebra structure in the
physical state space of an abstract TCFT (as well as in a
topological ``massive'' quantum field theory). Getzler's
ongoing work with J. Jones was already dealing with G-algebras
as the $n=2$ case of $n$-algebras.  Segal's ideas, see \cite
{se:cam}, played an essential role in Getzler's discovery. In particular,
Segal had already developed an extension to TCFTs of his geometric category
approach \cite{se:old} to conformal field theory.

        Later on Y.-Z. Huang \cite{h} found a third approach which
combined the ideas of Lian, Zuckerman and Getzler.  In
particular, Huang took steps towards the construction of a TCFT
from a TVOA; this work was based on Huang's earlier
demonstration of how to construct a (tree-level) holomorphic
conformal field theory (CFT) from a vertex operator algebra
(VOA). Such a connection between Segal's geometric approach
to conformal field theory and Borcherds' algebraic definition \cite{bor} of
a VOA (see also the book of I. Frenkel, J. Lepowsky and A.
Meurman \cite{flm}) had already been suggested in some public lectures by
I. Frenkel \cite{igor}.

        From the very beginning of the above development, it was
understood that the physical state space of a TCFT is merely the
cohomology of a much more enormous object, the BRST complex of the
TCFT.  Thus there arose the question: does the G-algebra structure on
the vector space of physical states come from a higher homotopy
G-algebra structure on the BRST complex itself?  This question was
explicitly raised by Lian and Zuckerman in their work on TVOAs and
associated TCFTs.  They found that in the BRST complex of a TVOA, all
of the identities of a G-algebra fail to hold on the nose, but they
continue to hold up to homotopy.  They then asked whether these
homotopies could be continued to an infinite hierarchy of higher
homotopies, such as those found in the work of Stasheff on
$A_\infty$-algebras and the more recent work of Stasheff and T. Lada on
$L_\infty$-algebras.

        The inspiration for the search for higher homotopy
algebras in topological conformal field theory arose in the
related context of closed string field theory, see Stasheff
\cite{jim:rec}. The explicit discussion of higher homotopy algebras
in string field theory appeared in work of Stasheff
\cite{jim:higher}, M. Kontsevich \cite{kon:sympl}, and E. Witten and
B. Zwiebach \cite{wz,z}. The later joint papers of T. Kimura,
A. Voronov and Stasheff \cite{ksv1,ksv2} constructed $L_\infty$
and $C_\infty$ structures using the operadic approach.  These
papers also include a conceptual explanation of the
relationship between string field theory and TCFTs.

        However, research on higher homotopies suffered from a lack of a
proper definition of a higher homotopy Gerstenhaber algebra (homotopy
G-algebra).  Recently, various definitions have been put forward, in
particular in work of V. Ginzburg and M. Kapranov \cite{gk}, Getzler and Jones
\cite{gj}, and Gerstenhaber and Voronov \cite{gv1}. In the
current paper, we use the term $G_\infty$-algebra to refer to a particular
definition of a \hg algebra appearing in the work of Getzler and Jones
(Definition~\ref{hg}). This definition is based on pioneering work of R.
Fox and L. Neuwirth. $G_\infty$-algebras are governed by what we call the
$G_\infty$-operad.

The main new result of the current paper is the proof that the
BRST complex of a TCFT is indeed a $G_\infty$-algebra.  In
fact,  a (tree-level) TCFT itself is defined to be an algebra
over a particular topological operad. Thus, in this paper, both
the ``classical''  theory of topological operads as well as the
recent theory of linear operads play essential roles.

        The original question of Lian and Zuckerman can now be
formulated precisely (see our Conjecture~\ref{conj}): does a TVOA
carry a natural $G_\infty$ structure? Since we have answered the
analogous question for a TCFT, a crucial step still remains in
the program to answer the original question: the completion of
Huang's work on the construction of a TCFT from a TVOA. Such a
construction should identify the BRST complex of the TCFT with
an appropriate topological completion of the BRST complex of
the TVOA. We look forward to the successful conclusion of this
program.

One of the essential tools of our paper is M. Wolf and
Zwiebach's ``string vertices'', which make a bridge between a
topological  operad of punctured Riemann spheres and the
infinite dimensional topological operad responsible for CFTs.
Amazingly, the latter operad plays a key role in the subjects
of Vassiliev invariants and double loop spaces. In particular,
string vertices combined with the approach of our paper yield
Vassiliev invariants of knots in Section~\ref{vass} and the
structure of a \hg algebra on the singular chain complex
of a double loop space in Section~\ref{double}.

\begin{sloppypar}
\begin{ack}
We are very grateful to T.~Q.~T.~Le, B.~Lian, A.~S. Schwarz, and J.~Stasheff
for helpful discussions. T.K. and A.A.V. express their sincere
gratitude to J.-L. Loday and J. Stasheff for their hospitality at the
wonderful conference in Luminy. T.K. and G.J.Z. would like to thank
A.A.V. for inviting them to the terrific conference at
Hartford. A.A.V. also thanks IHES for offering him excellent
conditions for work on the project in June of 1995.
\end{ack}
\end{sloppypar}

\section{Gerstenhaber algebras}
\label{g}

A \emph{Gerstenhaber algebra} or a \emph{G-algebra} is a graded vector
space $H$ with a dot product $xy$ defining the structure of a graded
commutative algebra and with a bracket $[x,y]$ of degree $-1$ defining
the structure of a graded Lie algebra, such that the bracket with an
element is a derivation of the dot product:
\[
[x , yz] = [x,y] z + (-1)^{(\deg x -1) \deg y} y [x,z] ,
\]
where $\deg x$ denotes the degree of an element $x$.  In other words,
a G-algebra is a specific graded version of a Poisson algebra.

This structure arises naturally in a number of contexts, such as the
following.

\begin{sloppypar}
\begin{ex}
Let $A$ be an associative algebra and $C^n(A, A) \linebreak[1] =
\linebreak[0]
\Hom ( A^{\otimes n} , A)$ be its Hochschild complex. Then
the dot product defined as the usual cup
product up to a sign
\begin{multline}
\label{dot-eq}
(x \cdot y) (a_1, \dots , a_{k+l}) = (-1)^{kl}
(x \cup y) (a_1, \dots , a_{k+l}) \\
= (-1)^{kl} x(a_1, \dots, a_k) y (a_{k+1}, \dots , a_{k+l}),
\end{multline}
where $x$ and $y$ are $k$- and $l$-cochains and $a_i \in A$, and a
G-bracket $[x,y]$ define the structure of a G-algebra on the \h
cohomology $H^n(A,A)$. The bracket was introduced by Gerstenhaber
\cite{gerst} in order to describe the obstruction for extending a
first order
deformation of the algebra $A$ to the second order. The following
definition of the bracket is due to Stasheff \cite{jim}. Considering
the tensor coalgebra $T(A) = \bigoplus_{n=0}^\infty A^{\otimes n}$
with the comultiplication $\Delta (a_1 \otimes \dots \otimes a_n) =
\sum_{k=0}^n (a_1 \otimes \dots \otimes a_k) \otimes (a_{k+1} \otimes
\dots \otimes a_n)$, we can identify the Hochschild cochains $\Hom
(A^{\otimes n}, A)$ with the coderivations $\coder T(A)$ of the tensor coalgebra
$T(A)$. Then the G-bracket $[x,y]$ is defined as the (graded)
commutator of coderivations. In fact, the \h complex $C^\bullet
(A,A)$ is a differential graded Lie algebra with respect to
this bracket.
\end{ex}
\end{sloppypar}

\begin{ex}
Let $A^\bullet_n$ be the $\nz$-graded commutative algebra generated by $n$
variables $x_1, x_2, ..., x_n$, of degree zero, and $n$ more variables
$\partial_{x_1}, \dots, \partial_{x_n}$ of degree one. We refer to an element
of this algebra as a polyvector field. The elements of degree zero are
interpreted as functions, the elements of degree one as vector fields, those
of degree two as bivector fields, and so on.  The dot product is simply the
graded commutative multiplication of polyvector fields.

        Long ago, Schouten and Nijenhuis \cite{nijenhuis} defined a bracket
operation on polyvector fields (they thought of such fields as
antisymmetric contravariant tensor fields). The Schouten-Nijenhuis
bracket $[P,Q]$ is characterized by the
following:

\begin{enumerate}
\item For any two functions $f$ and $g$, $[f,g] = 0$.

\item  If $f$ is a function and $X$ is a vector field, $[X,f] = -[f,X] = Xf$.

\item If $X$ and $Y$ are vector fields, then $[X,Y]$ is the standard
bracket of the vector fields.

\item Together, the dot product and the Schouten-Nijenhuis bracket endow
$A^\bullet_n$ with the structure of a Gerstenhaber algebra.
\end{enumerate}

Let $C_n$ be the polynomial algebra in $n$ variables.  It is known that
$H^\bullet (C_n, C_n)$ is canonically isomorphic as a G-algebra to the
algebra $A^\bullet_n$.
\end{ex}

\begin{ex}
\label{3}
 Let $\gtg$ be any Lie algebra, and let $\Lambda^\bullet \gtg$
be the Grassmann algebra generated by $\gtg$.  Define a bracket $[X,Y]$ on
$\Lambda^\bullet \gtg$ by requiring the following:

\begin{enumerate}
\item If $a$ and $b$ are scalars, $[a,b] = 0$.

\item If $X$ is in $\gtg$ and $a$ is a scalar, then $[X,a] = 0$.

\item If $X$ and $Y$ are in $\gtg$, then  $[X,Y]$ is the Lie bracket in $\gtg$.

\item The wedge product together with the bracket product endow
$\Lambda^\bullet \gtg$ with the structure of a Gerstenhaber algebra.
\end{enumerate}
\end{ex}

\begin{ex}
Let $M$ be a manifold (differentiable, complex, algebraic,
etc.)  Let $F(M)$ be the commutative algebra of functions (of the
appropriate type--differentiable, holomorphic, regular, etc.) on $M$.
Let $V^\bullet(M)$ be the algebra of polyvector fields on $M$, with the
operations of wedge product and the Schouten-Nijenhuis bracket, defined by
analogy with the bracket in $A^\bullet_n$.  Then $V^\bullet(M)$ is a
Gerstenhaber algebra.  We can regard $V^\bullet(M)$ as the commutative
superalgebra of functions on $\Pi T^* M$, the cotangent bundle of $M$
with the fibers made into odd supervector spaces.  The G-bracket in
$V^\bullet(M)$ is the odd Poisson bracket associated to the canonical
odd symplectic two-form on $\Pi T^* M$.

         Let $P$ be a bivector field on $M$.  We can always construct a
bracket operation on the algebra $F(M)$ by the formula
\[
        \{f,g\}_P = \iota(P) (df \wedge dg) =  (df \wedge dg) (P),
\]
where $\iota(P)$ denotes contraction of $P$ against a two-form.  The bracket
$\{ , \}_P$ satisfies the Jacobi identity if and only if the
Schouten-Nijenhuis
bracket $[P,P]$ is zero.  In this case, the algebra $F(M)$ becomes what is
known as a Poisson algebra; moreover, the derivation $\sigma_P =
[P,-]$ has square zero and turns $V^\bullet (M)$ into a differential
graded G-algebra, whose cohomology G-algebra is known as the Poisson
cohomology of $M$ relative to $P$.
\end{ex}

\begin{ex}
Let $X$ be a topological space and let $\Omega^2 X$
be the two-fold loop space of $X$.  Let $A^\bullet (X)$ denote
the homology of $\Omega^2 X$ with rational coefficients.  We
endow $A^\bullet (X)$ with the structure of a $\nz$-graded
commutative algebra, via the Pontrjagin product. $A^\bullet
(X)$ is a Hopf algebra which is freely generated as a graded
commutative algebra by its subspace $P$ of primitive elements.
Moreover, $P$ is isomorphic to the rational homotopy of
$\Omega^2 X$.  $P$ is $\nz$-graded, and we have $P_n = \pi_{n+2}(X)
\otimes_\nz \nq$, for $n$ nonnegative, and 0 otherwise.

The rational homotopy of the ordinary loop space of $X$ is a $\nz$-graded Lie
algebra, which we denote by $L$.  The bracket is called the Samelson product.
It is known that $P_n = L_{n+1}$.  Thus, the algebra $A^\bullet (X)$ is
isomorphic to the graded exterior algebra $\Lambda^\bullet L$.  In particular,
$A^\bullet (X)$ is a G-algebra of a type generalizing Example~\ref{3} above.
\end{ex}

\subsection{Operads in action}

For a primer on operads, algebras over operads and the little
disks operad, see P. May's paper \cite{may:def} in this volume.
Here we recall briefly the definition of the little disks
operad in relation to G-algebras. The \emph{little disks operad}
is the collection $\{D(n), \, n \ge 1\}$ of topological spaces
$D(n)$ with an action of the permutation group and operad
compositions. The space $D(n)$ consists of embeddings of $n$
little disks in the unit disk via dilatation and translation.
The operad composition $\circ_i: D(m) \times D(n) \to D(m+n-1)$
is defined by contracting the unit disk with $m$ little disks
inside it to fit into the $i$th little disk in the other unit
disk and erasing the seam. Since $D(n)$ is naturally an open
subset in $\nr^{3n}$, it is a topological operad and one can
naturally obtain an operad of graded vector spaces from it,
applying the functor of homology. This operad $H_\bullet
(D(n))$ is called a \emph{G-operad} in view of the following
theorem.

\begin{thm}[F. Cohen \cite{C1,fc}]
The structure of a G-algebra on a $\nz$-graded vector space is
equivalent to the structure of an algebra over the homology
little disks operad $H_\bullet(D(n))$.
\end{thm}

\section{Topological vertex operator algebras}
\label{tvoa}

In this section, we give a brief introduction to VOAs, in order
to present the Lian-Zuckerman conjecture in an updated form.
VOAs will not show up in the rest of the paper, but we intend
to get back to them in a subsequent paper. See \cite{lz} for more discussion.

\begin{sloppypar}
\begin{df}[Quantum operators]
Let $V^\bullet [\cdot] = \bigoplus_{g \in \nz, \, \Delta \in \nz}
V^g[\Delta]$ be an integrally bigraded complex vector space; if $v$ is in
$V^g[\Delta]$ we will write $|v| = g = \text{the ghost number of }v$ and
$||v|| = \Delta = \text{the weight of }v$.  Let $z$ be a formal variable with
degrees $|z|= 0$ and $||z|| = -1$.  Then, it makes sense to speak of a
homogeneous \emph{bi-infinite} formal power series
\[
\phi(z) = \sum_{n \in \nz} \phi(n) z^{-n-1}
\]
of degrees $|\phi(z)|$, $||\phi(z)||$, where the coefficients $\phi(n)$ are
homogeneous linear maps in $V^\bullet [\cdot]$ of degrees $|\phi(n)| =
|\phi(z)|$, $||\phi(n)|| = -n-1+||\phi(z)||$.  Note then that the terms
$\phi(n) z^{-n-1}$ indeed have the same degrees $|\phi(z)|$, $||\phi(z)||$
for all $n$.  We call a finite sum of such series a \emph{quantum operator}
on $V^\bullet [\cdot]$, and denote the bigraded linear space of quantum
operators as $\QO(V^\bullet [\cdot])$.  We will denote the special operator
$\phi(0)$ by the symbol $\Res_z \phi (z)$. 
\end{df}
\end{sloppypar}

\begin{df}
A \emph{vertex operator graded algebra} consists of the
following ingredients:
\begin{enumerate}
\item An integrally bigraded complex vector space $V^\bullet [\cdot]$.

\item A linear map $Y: V^\bullet [\cdot] \to \QO(V^\bullet [\cdot])$
such that $Y$ has bidegree (0,0).  We call $Y$ the \emph{vertex map}, and if
$v$ is
in $V^\bullet [\cdot]$, we let $Y(v,z) = Y(v) (z)$ denote the \emph{vertex
operator} associated to $v$. The map $Y$ is subject to the following axioms:
   \begin{enumerate}
   \item Let $v$ and $v'$ be elements of $V^\bullet [\cdot]$: then $\Res_z
   (z^m Y(v,z)v')$ vanishes for $m$ sufficiently positive.

   \item (Cauchy-Jacobi identity) Let $v$ and $v'$ be elements of $V^\bullet
   [\cdot]$ and let $f(z, w)$ be a Laurent polynomial in $z$, $w$, and $z-w$.
   Then we have the identity
\noindent
  \begin{multline*}
   \Res_w \Res_{z-w} Y(Y(v, z-w)v', w)f(z,w) \\
   = \Res_z \Res_{w} Y(v,z) Y(v',w) f(z,w) \\
   - (-1)^{|v|\,|v'|} \Res_w \Res_{z} Y(v',w) Y(v, z) f(z,w) .
   \end{multline*}
   \item There exists a distinguished element 1 in $V^0[0]$ such that $Y(1, z)$
   is the
   identity operator and such that for any $v$ in $V^\bullet [\cdot]$,
   the result $Y(v, z)1$ is a power series in $z$ and
   \[
   \lim_{z \to 0} Y(v,z) 1 = v .
   \]
   \end{enumerate}

\item A distinguished element $F$ in $V^0[1]$ such that $F_0 = \Res_z Y(F,
z)$ defines the ghost number grading: if $v$ is in $V^g[\Delta]$, $F_0
v= gv$.

\item A distinguished element $L$ in $V^0[2]$ such that if we define an
operator $L_n = \Res_z (z^{n+1} Y(L,z))$ for every integer $n$,
we have the following:
\begin{enumerate}
\item  For $v$ in $V^g[\Delta]$, $L_0 v = \Delta v$.

\item For any $v$ we have $Y(L_{-1} v, z) = \partial Y(v,z)$.

\item For some fixed complex number $c$ (the central charge), we have
\[
[L_m, L_n] = (m-n) L_{m+n}  + \frac{c}{12} (m^3 - m) \delta_{m, -n} 1 .
\]
\end{enumerate}
\end{enumerate}
\end{df}

\begin{df}
A \emph{topological vertex operator algebra $($TVOA}) is a
vertex operator graded algebra $(V^\bullet [\cdot], Y, 1, F, L)$
equipped with two additional distinguished elements $J$ in $V^1[1]$ and
$G$ in $V^{-1} [2]$ such that the following axioms hold:

  \begin{enumerate}
  \item The operator $Q = \Res_z Y(J,z)$ satisfies $Q^2 = 0$. $Q$ is called
the BRST charge, or BRST coboundary operator.  $Q$ has bidegree
$(1, 0)$: $|Q| = 1$, $||Q|| = 0$.

  \item $[Q, Y(G, z)] = Y(L, z)$.
  \end{enumerate}
\end{df}

\begin{df}
Let $(V^\bullet [\cdot], Y, 1, F, L, J, G)$ be a
TVOA.  Let $v$ and $v'$ be elements of $V^\bullet [\cdot]$.

  \begin{enumerate}
  \item The \emph{dot product}  of $v$ and $v'$ is the element
\[
v \cdot v' = \Res_z (z^{-1} Y(v,z)v').
\]

  \item The \emph{bracket product}  of $v$ and $v'$ is the element
\[
        [v, v'] =(-1)^{|v|} \Res_w \Res_{z-w} Y(Y(G,z-w)v, w) v' .
\]
  \end{enumerate}
\end{df}

\begin{lm}
Let $V^\bullet [\cdot]$ be a TVOA.
\begin{enumerate}
\item The BRST operator $Q$ is a derivation of both the dot and bracket
products, which therefore induce dot and bracket products respectively
on the BRST cohomology $H^\bullet (V, Q)$ of $V^\bullet [\cdot]$
relative to $Q$.

\item Every BRST cohomology class is represented by a BRST-cocycle of
weight 0.  Thus, the weight grading induces the trivial grading in
BRST cohomology.

\item The BRST cohomology is trivial unless the central charge $c = 0$.
\end{enumerate}
\end{lm}

\begin{thm}
Let $V^\bullet [\cdot]$ be a TVOA.  With respect to to the
induced dot and bracket products, the BRST cohomology $H^\bullet (V,
Q)$ of $V^\bullet [\cdot]$ relative to $Q$ is a Gerstenhaber algebra.
\end{thm}

For an interesting study of the identities satisfied by the dot
and bracket products in the TVOA $V^\bullet [\cdot]$ itself,
see F.~Akman \cite{akman}. In particular, Akman finds that the
space $V^\bullet [\cdot]$ endowed with the bracket is a
$\nz$-graded Leibniz algebra \cite{loday}.

\begin{conj}
\label{conj}
Let $V^\bullet [\cdot]$ be a TVOA.  Then the dot product and the
skew-symmetrization of the bracket defined above can be extended to
the structure of a $G_\infty$-algebra (see below) on $V^\bullet
[\cdot]$.
\end{conj}

This conjecture makes precise sense in the light of our current
work: it refines the question posed by Lian and Zuckerman, who
expected $A_\infty$ and $L_\infty$ structures to mix together.

\section{Classical story: homotopy associative and homotopy Lie
algebras}

Inasmuch as G-algebras combine properties of commutative associative
and Lie algebras, \hg algebras make a similar combination of homotopy
associative and \hl algebras; see Section~\ref{hg-op} below. (The
commutativity is not completely lost, either: \hg algebras will also provide a
homotopy for the commutativity of the dot product). Before
discussing \hg algebras, let us recall definitions of the more
traditional homotopy associative and \hl algebras.

\begin{df}[Homotopy associative ($A_\infty$-) algebras]
A {\it homotopy associative algebra} is a complex $V = \sum_{i \in
\nz} V_i$ with a differential $d$, $d^2 = 0$, of degree 1 and a
collection of $n$-ary products $M_n$:
\[
M_n(v_1, \dots, v_n) \in V, \qquad v_1, \dots, v_n \in V,\; n \ge 2,
\]
which are homogeneous of degree $2-n$ and satisfy the relations
\begin{multline*}
d M_n(v_1, \dots, v_n) + \sum_{i=1}^n  \epsilon(i) M_n(v_1, \dots, dv_i, \dots,
v_n)
\\
= \sum_{
  \substack{
   k+l = n+1\\
   k, l \ge 2
}}
\sum_{i=0}^{l-1}
\epsilon(k,i)
M_l(v_1, \dots, v_i, M_k(v_{i+1}, \dots, v_{i+k}), v_{i+k+1}, \dots, v_{n}),
\end{multline*}
where $\epsilon(i) = (-1)^{\deg v_1 + \dots + \deg v_{i-1}}$ is the
sign picked up by taking $d$ through $v_1, \linebreak[0] \dots,
\linebreak[1] v_{i-1}$, $\epsilon(k,i) = (-1)^{k(\deg v_1 + \dots +
\deg v_{i})}$ is the sign picked up by $M_k$ passing through $v_1,
\linebreak[0] \dots, \linebreak[1] v_{i}$.
\end{df}
For $n=3$, the above identity shows that the binary product $M_2$ is
associative up to a homotopy, provided by the ternary product $M_3$.

\begin{df}[Homotopy Lie ($L_\infty$-) algebras]
A {\it homotopy Lie algebra} is a complex $V = \sum_{i \in \nz} V_i$ with a
differential $d$, $d^2 = 0$, of degree 1 and a collection of $n$-ary brackets:
\[
[v_1, \dots, v_n] \in V, \qquad v_1, \dots, v_n \in V,\; n \ge 2,
\]
which are homogeneous of degree $3-2n$ and super (or graded) symmetric:
\[
[v_1, \dots, v_i,  v_{i+1}, \dots , v_n] = (-1)^{|v_i| |v_{i+1}|} [v_1, \dots,
v_{i+1}, v_i, \dots , v_n],
\]
$\deg v$ denoting the degree of $v \in V$,
and satisfy the relations
\begin{multline*}
d[v_1, \dots, v_n] + \sum_{i=1}^n  \epsilon(i) [v_1, \dots, dv_i, \dots, v_n]
\\
= \sum_{
  \substack{
   k+l = n+1\\
   k, l \ge 2
}}
\sum_{
  \substack{
   \text{unshuffles } \sigma:\\
   \{1,2, \dots, n\} = I_1 \cup I_2,\\
   I_1 = \{i_1, \dots, i_k\}, \; I_2 = \{ j_1, \dots, j_{l -1}\}
}}
\epsilon (\sigma)
[[v_{i_1}, \dots, v_{i_k}], v_{j_1}, \dots, v_{j_{l-1}}],
\end{multline*}
where $\epsilon (i) = (-1)^{\deg v_1 + \dots + \deg v_{i-1}}$ is the
sign picked up by taking $d$ through $v_1, \linebreak[0] \dots,
\linebreak[1] v_{i-1}$, $\epsilon (\sigma)$ is the sign picked up by
the elements $v_i$ passing through the $v_j$'s during the unshuffle of
$v_1, \dots , v_n$, as usual in superalgebra.
\end{df}
For $n=3$, the above identity shows that the binary bracket $[v_1,
v_2]$ satisfies the Jacobi identity up to a homotopy, provided by the
next bracket $[v_1,v_2,v_3]$.

\section{Homotopy G-algebras}

Due to Fred Cohen's Theorem \cite{fc}, the structure of a G-algebra on
a vector space $V$ is equivalent to the structure on $V$ of an algebra
over the homology operad $H_\bullet (D(n))$, $n \ge 1$, of the little
disks operad $D(n)$. In most general terms, a homotopy Gerstenhaber
algebra or homotopy G-algebra is an algebra over a ``resolution'' of
this G-operad $G_\bullet (n) = H_\bullet (D(n))$, $n \ge 1$, i.e., an
operad $hG_\bullet(n)$ of complexes whose cohomology is identified
with $G_\bullet (n)$. Different resolutions lead to different notions
of homotopy G-algebras. A minimal resolution to answer the question of
Deligne that the \h complex is a \hg algebra was
introduced in \cite{gv1}, two free resolutions were constructed in
\cite{gj}.  The obvious singular-chain resolution $C_\bullet (D(n))$
is too large for our purposes: it produces too many higher operations
(homotopies).

\begin{ex}[Hochschild complex of an associative algebra]
\label{hochschild}
Let $A$ \linebreak[4]
be an associative algebra and $C^n (A, A) = \Hom
(A^{\otimes n}, A)$ its \hc\ complex.
Define the following
collection of multilinear operations, called {\it braces,} on $C^\bullet(A,A)$:
\begin{multline*}
\{x\} \{ x_1, \dots, x_n\} (a_1, \dots, a_m)
 := \\
  \sum
(-1)^\varepsilon x(a_1, \dots , a_{i_1},  x_1 (a_{i_1+1}, \dots), \dots,
a_{i_n},
  x_n(a_{i_n+1}, \dots),  \dots, a_m)
\end{multline*}
for $x , x_1, \dots , x_n \in C^\bullet(A,A)$, $a_1, \dots, a_m \in
A$, where the summation runs over all possible substitutions of $x_1,
\dots, x_n$ into $x$ in the prescribed order and $\varepsilon :=
\sum_{p=1}^{n} (\deg x_{p}-1) i_{p} $. The braces $\{x\}\{ x_1, \dots,
x_n\}$ are homogeneous of degree $-n$ , i.e., $\deg \{x\}\{ x_1, \dots,
x_n\} = \deg x + \deg x_1 + \dots + \deg x_n -n $. We will also adopt
the following convention:
\[
x \circ y : = \{x\} \{y\}.
\]
In addition, the usual cup product (altered by a sign) defined by
\eqref{dot-eq} and the differential
\begin{multline*}
(d x) (a_1, \dots, a_{n+1}) \\
\begin{split}
 := \; &
\;(-1)^{\deg x} a_1 x(a_2, \dots, a_{n+1}) \\
 &  + (-1)^{\deg x}  \sum_{i=1}^n (-1)^i x (a_1, \dots , a_{i-1}, a_i a_{i+1},
a_{i+2}, \dots , a_{n+1})\\
 & - x(a_1, \dots, a_n) a_{n+1}
\end{split}
\end{multline*}
define the structure of a differential graded (DG) associative algebra
on $C^\bullet (A,A)$. In their turn the braces satisfy the following
identities.
\begin{multline}
\label{higher}
\{\{x\} \{x_1, \dots, x_m\}\} \{ y_{1}, \dots , y_{n}\} \\
=  \sum_{0 \le i_1 \le \dots \le i_m \le n}  (-1)^\varepsilon  \{x\} \{ y_1,
\dots,
y_{i_1}, \{x_1\}  \{ y_{i_1+1}, \dots \},  \dots, \\
 y_{i_{m}}, \{x_m \} \{  y_{i_{m}+1} ,
\dots \}, \dots, y_n \},
\end{multline}
where $\varepsilon := \sum_{p=1}^{m} (\deg x_{p} - 1) \sum_{q=1}^{i_p}
(\deg y_{q} - 1)$, i.e., the sign is picked up by the $\{x_i\} $'s
passing through the $\{y_j\}$'s in the shuffle.
\begin{equation}
\label{dist}
\{x_1 \cdot x_2\} \{ y_1, \dots, y_n\} = \sum_{k=0}^n (-1)^\varepsilon \{x_1\}
\{ y_1, \dots, y_k\} \cdot \{x_2\} \{ y_{k+1}, \dots, y_n\},
\end{equation}
where $\varepsilon = (\deg x_2) \sum_{p=1}^k (\deg y_p - 1 )$.
\begin{multline}
\label{comm}
\begin{split}
&d ( \{x\} \{ x_1, \dots, x_{n+1}\} )
- \{dx\}\{x_1, \dots, x_{n+1}\} \\
&- (-1)^{\deg x -1} \sum_{i=1}^{n+1} (-1)^{\deg x_1 + \dots + \deg
x_{i-1} - i - 1 }   \{x\} \{ x_1,
\dots, dx_i, \dots , x_{n+1}\}
\end{split}
\\
\begin{split}
\; = \; & - (-1)^{\deg x (\deg x_1 -1) } x_1 \cdot \{x\} \{ x_2, \dots,
x_{n+1} \} \\ & - (-1)^{\deg x} \sum_{i=1}^n (-1)^{\deg x_1 + \dots +
\deg x_{i} - i } \{x\} \{ x_1, \dots , x_i \cdot x_{i+1}, \dots , x_{n+1}
\}\\ & + (-1)^{\deg x + \deg x_1 + \dots + \deg x_{n} - n} \{x\} \{ x_1,
\dots, x_n \} \cdot x_{n+1}
\end{split}
\end{multline}
The structure of braces $ \{x\} \{x_1, \dots, x_n\} $, $n \ge 1$, and a
dot product $xy$ satisfying the above identities on a complex $V =
\bigoplus_n V^n$ of vector spaces was called in \cite{gv1} a {\it \hg
algebra\/}. Thus the above braces and dot product provide the
Hochschild complex $C^\bullet (A, A)$ of an associative algebra $A$
with with the structure of a \hg algebra.  From the definition of a
\hg algebra, it is easy to find the definition of a \hg operad: it is
an operad of complexes with braces $\{x\}\{x_1, \dots, x_n\}$ and a dot
product $xy$ as generators and the associativity and the above
identities for the braces and the dot product as defining relations.
\end{ex}

\begin{rem}
\label{B-infty}
A more general notion of a \hg algebra is that of $B_\infty$-algebra; 
see Getzler and Jones \cite{gj}. A $\nz$-graded vector space
$V^\bullet$ is a $B_\infty$-algebra if and only if the bar coalgebra
$\BV = \bigoplus_{n=0}^\infty (V[-1])^{\otimes n},$ has the structure
of a DG bialgebra, that is to say, a bialgebra with a differential
which is simultaneously a derivation and a coderivation. The
differential gives rise to a differential on $V$, as well as the dot
product on $V$ and higher $A_\infty$ products. The multiplication,
compatible with the coproduct, on $\BV$ gives rise to braces on $V$ and
some higher braces. If the higher products and the higher braces
vanish identically, we retrieve the \hg algebra structure of the
previous example.
\end{rem}

In this paper we will need the following still more general geometric
definition of a \hg algebra of Getzler and Jones \cite{gj}; see also
\cite{gv2}.  Consider Fox-Neuwirth's cellular partition of the
configuration spaces $F(n,
\nc) = \nc^n \setminus \Delta$, where $\Delta$ is the weak diagonal
$\cup_{i \ne j} \{x_i = x_j\}$, of $n$ distinct points on the complex
plane $\nc$: cells are labeled by ordered partitions of the set $\{1,
\dots, n\}$ into subsets with arbitrary reorderings within each subset.
This
reflects grouping points lying on common vertical lines on the plane
and ordering the points lexicographically. For each $n$, take the
quotient cell complex $K_\bullet \mathcal M (n)$ by the action of
translations $\nr^2$ and dilations $\nr^*_+$. These quotient spaces do
not form an operad, but one can glue lower $K_\bullet \mathcal M
(n)$'s to the boundaries of higher $K_\bullet \mathcal M (n)$'s to
form a cellular operad $K_\bullet \FF = \{K_\bullet \F{n}\; | \; n \ge
2\}$. In fact, the underlying spaces $\F{n}$ are manifolds with corners
compactifying $\mathcal{M}(n)$.

The resulting space $\F{n}$ fibers over the real
compactification $\X{n}$ of the moduli space $\mathcal M_{0,n}$ of
$n$-punctured curves of genus zero, see \cite{gv2,gj}. The
space $\F{n}$ can be also interpreted as a ``decorated'' moduli space, see
\cite{gv2}.
Indeed, it can be identified with the moduli space of data $(C; p_1, \dots,
p_{n+1}; \linebreak[0] \tau_1, \dots, \tau_m, \tau_\infty)$, where $C$
is a stable algebraic curve with $n+1$ punctures $p_1, \dots, p_{n+1}$
and $m$ double points. For each $i$, $ 1 \le i \le m$, $\tau_i$ is the
choice of a tangent direction at the $i$th double point to the
irreducible component that is furthest away from the ``root'', i.e.,
from the component of $C$ containing the puncture $\infty := p_{n+1}$,
while $\tau_\infty$ is a tangent direction at $\infty$. The stability of a
curve is understood in the sense of Mumford's geometric invariant
theory: each irreducible component of $C$ must be stable, i.e., admit
no infinitesimal automorphisms. The operad composition is given by
attaching the $\infty$ puncture on a curve to one of the other punctures on
another curve, remembering the tangent direction at each new double
point.

Cells in this cellular operad $K_\bullet \FF$ are
enumerated by pairs $(T, p)$, where $T$ is a directed rooted tree with
$n$ labeled initial vertices and one terminal vertex, labeling a
component of the boundary of $\F{n}$, and $p$ is an ordered partition,
as above, of the set $\IN (v)$ of incoming edges for each vertex
$v$ of the tree $ T$.

In \cite{gj}, it is shown that a complex $V$ is a \hg algebra of the
example above iff it is
an algebra over the operad $K_\bullet \FF$ satisfying the following condition.
The structure mappings
\[
K_\bullet \F{n} \to \Hom (V^{\otimes n}, V),
\]
for the algebra $V$ over the operad $K_\bullet \FF$ send all cells in
$K_\bullet 
\mathcal M (n)$ to zero, except cells of two kinds:
\begin{enumerate}
\item $(\delta_n; i_1| i_2 \dots i_n)$, where $\delta_n$ is the corolla,
the tree with one root and $n$ edges, connecting it to the remaining $n$
vertices, corresponding to the configuration where the points $i_2, \dots, i_n$
sit on a vertical line, the $i_k$th point being below the $i_{k+1}$st, and the
$i_1$st point is in the half-plane to the left of the line;
\item $(\delta_2; i_1 i_2)$,
corresponding to the configuration where all the points sit on a
single vertical line, the $i_1$st point being below the $i_{2}$nd.
\end{enumerate}
Cells of the first kind give rise to the braces $\{x_{i_1}\} \{ x_{i_2},
\dots , x_{i_n} \}$, $n \ge 2$, and cells of the second kind give rise
to the dot products $x_{i_1} x_{i_2}$.  The relations \eqref{dist} and
\eqref{comm}, and the associativity of the dot product, follow from
the
combinatorial structure of the cell complex. The relation
\eqref{higher} does not necessarily follow from the cell complex. Thus
the \hg operad of Example~\ref{hochschild} is a quotient operad of
$K_\bullet \FF$. In particular, a \hg algebra of Example~\ref{higher}
is an algebra over $K_\bullet \FF$. In fact, the $B_\infty$-operad
corresponding to $B_\infty$-algebras of Remark \ref{B-infty} is an
intermediate operad between $K_\bullet \FF$ and the \hg operad of
Example~\ref{hochschild}. Thus, every $B_\infty$-algebra is a $K_\bullet
\FF$-algebra
with certain operations also set to zero, and every \hg algebra of
Example~\ref{hochschild} is a $B_\infty$-algebra with some operations
vanishing.

The main purpose of this
paper is to give a nontrivial example of an algebra over the operad
$K_\bullet
\FF$ and thus answer the question of Lian and Zuckerman \cite{lz}
about the homotopy structure of a TCFT.
We will adopt the following definition of a \hg operad and a \hg
algebra, introduced in \cite{gj} under the names of a homotopy
2-algebra and a braid algebra. We think it is appropriate to
call them the $G_\infty$-operad and a $G_\infty$-algebra,
respectively.
\begin{df}
\label{hg}
The \emph{$G_\infty$-operad} is the operad $K_\bullet \FF$ of complexes. An
algebra over it is called a \emph{$G_\infty$-algebra}.
\end{df}

\subsection{Some lower operations}

Let $V$ be a $G_\infty$-algebra, as in Definition~\ref{hg}. This structure
assumes a collection of $n$-ary operations $\mu_{T,p}$ on $V$ associated with
each cell of $K_\bullet \F{n}$, i.e., with each pair $(T,p)$, where $T$ is an
$n$-tree labeling a component of the boundary of $\F{n}$ and $p$ is an ordered
partition of incoming vertices for each vertex of $T$, see above. Here we would
like to illustrate this complicated algebraic  structure with an explicit
description of the operations and relations for $n =2$ and 3.

For $n=2$, $\F{2} = \MM (2) \cong S^1$ and there are just two types of
cells: $(\delta_2; 1|2)$ and $(\delta_2; 12)$, dividing the
circle into two intervals and two points. Denote the corresponding
operations on $V$ by $v_1 \circ v_2 = \{v_1\}\{v_2\}$ and $v_1 v_2$,
where $v_1$ and $v_2 \in V$. The cells $(\delta_2; 2|1)$ and
$(\delta_2; 21)$ are obtained by applying the transposition
$\tau_{12} \in S_2$ to the above two, and therefore the corresponding
operations are just $v_2 \circ v_1$ and $v_2 v_1$.

We need also to introduce the following bracket, which, unlike
the circle product, induces an operation on the $d$-cohomology:
\begin{equation}
\label{bracket}
[v_1,v_2] = v_1 \circ v_2 - (-1)^{(\deg v_1 - 1)(\deg v_2 -1)} v_2 \circ v_1 .
\end{equation}

For $n=3$, $\MM (3) \subset \F{3}$ contains cells of the following types:
$(\delta_3; 1|2|3)$, $(\delta_3; 1|23)$, $(\delta_3; 12|3)$, and
$(\delta_3; 123)$, corresponding to operations which we will denote by
$\{v_1\}\{v_2\}\{v_3\}$, $\{v_1\}\{v_2,v_3\}$, $\{v_1,v_2\}\{v_3\}$, and
$M_3(v_1,v_2,v_3)$, respectively.

\subsection{Some lower identities}
\label{identities}

There are no identities between compositions of operations, because
the operad $\FF$ is free as an operad, and so is $K_\bullet \FF$.  But
there are identities involving the differential $d$, because the
boundary of a cell in $K_\bullet \F{n}$ is a linear combination of
other cells --- $K_\bullet \FF$ is a cellular operad after all. The
sign rule we use for the boundaries of cells in the sequel is as
follows. We introduce an orientation on cells in the configuration space
$\F{n}$, ordering their coordinates according to the rule: first,
going from left to right, the $x$ coordinates of the lines on which
the points group, then, going lexicographically from left to right
and from top to bottom, the $y$ coordinates
of the points in each group.

The boundary of the 1-cell $(\delta_2; 1|2) \in K_\bullet \F{2}$
is

\let\picnaturalsize=N
\def\picsize{3.0in}
\def\picfilename{d2.eps}
\ifx\nopictures Y\else{\ifx\epsfloaded Y\else\input epsf \fi
\let\epsfloaded=Y
\centerline{\ifx\picnaturalsize N\epsfxsize \picsize\fi
\epsfbox{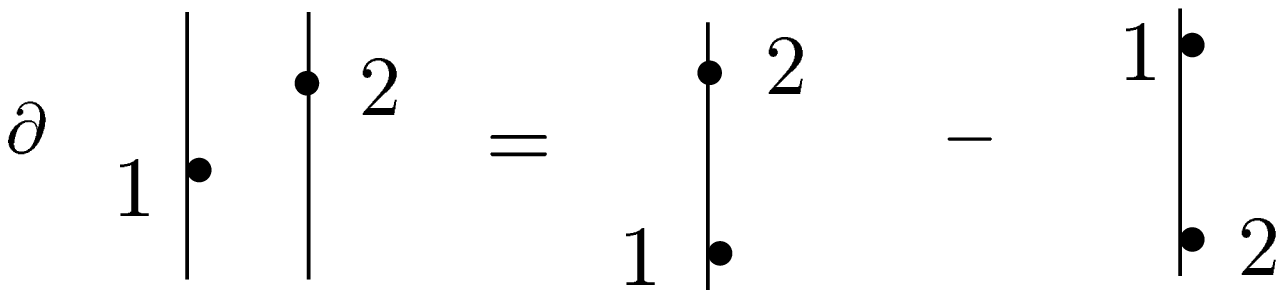}}}\fi
\noindent
\[
\partial (\delta_2; 1| 2)  =  (\delta_2; 1 2) - (\delta_2; 2 1).
\]
This incidence relation between cells implies the following \emph{homotopy
commutativity} relation for the dot product:
\begin{multline}
\label{hcomm}
d(v_1 \circ v_2) - dv_1 \circ v_2 - (-1)^{\deg v_1 - 1} v_1 \circ dv_2
\\= v_1 v_2 - (-1)^{(\deg v_1)(\deg v_2)} v_2 v_1.
\end{multline}
We also have $\partial (\delta_2; 12) = 0$, which yields the
following \emph{derivation property} of $d$ with respect to the dot
product:
\[
d(v_1 v_2) - dv_1  v_2 - (-1)^{\deg v_1} v_1 dv_2 = 0.
\]

For the top cell in $K_\bullet \F{3}$, we have

\let\picnaturalsize=N
\def\picsize{3.0in}
\def\picfilename{d3.eps}
\ifx\nopictures Y\else{\ifx\epsfloaded Y\else\input epsf \fi
\let\epsfloaded=Y
\centerline{\ifx\picnaturalsize N\epsfxsize \picsize\fi
\epsfbox{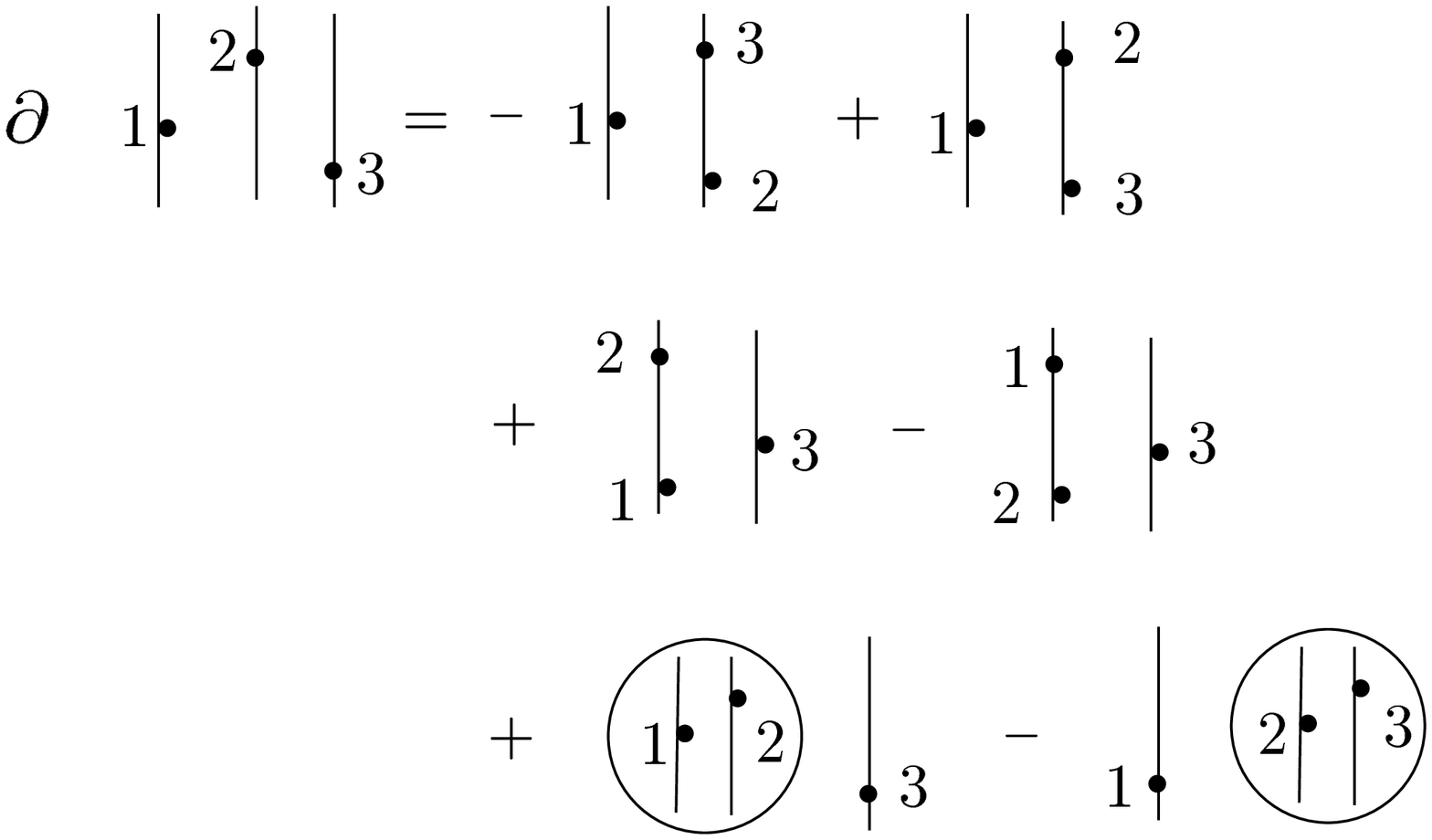}}}\fi
\noindent
\begin{align*}
\partial (\delta_3; 1| 2| 3)
 =  & -   (\delta_3; 1|23) + (\delta_3;
1| 3 2) \\
& + (\delta_3; 12|3)
- (\delta_3; 2 1| 3)\\
& + (\delta_2
\circ_1 \delta_2, (1 |  2) \circ_1 (1| 2)) - (\delta_2 \circ_2
\delta_2; (1|2) \circ_2 (1| 2)),
\end{align*}
the circles on the figure meaning ``magnifying glasses'', here showing
that a pair of points on the sphere bubbled off onto another
``microscopic'' sphere attached at a double point. This happens in
much the same way that a sphere pinches off to form a genus zero curve
with a node as one goes to the compactification divisor in the moduli
space of genus zero curves. This implies
\begin{multline}
\label{distr}
\begin{split}
d (\{v_1\}\{v_2\}\{v_3\}) - \{ d v_1\}\{v_2\}\{v_3\} - (-1)^{\deg v_1
- 1} \{v_1\}\{d v_2\}\{v_3\} \\
 - (-1)^{\deg v_1 + \deg v_2}
\{v_1\}\{v_2\}\{d v_3\}
\end{split}
\\
\begin{split}
=& - \{v_1\}\{v_2,v_3\} + (-1)^{(\deg v_2-1)(\deg v_3
-1)} \{v_1\}\{v_3,v_2\} \\
&+ \{v_1,v_2\}\{v_3\}
- (-1)^{(\deg v_1 -1)(\deg v_2
-1)} \{v_2,v_1\}\{v_3\} \\
& +
 ( v_1 \circ v_2) \circ v_3 - v_1 \circ (v_2
\circ v_3).
\end{split}
\end{multline}
The incidence relations

\let\picnaturalsize=N
\def\picsize{3.0in}
\def\picfilename{d31.eps}
\ifx\nopictures Y\else{\ifx\epsfloaded Y\else\input epsf \fi
\let\epsfloaded=Y
\centerline{\ifx\picnaturalsize N\epsfxsize \picsize\fi
\epsfbox{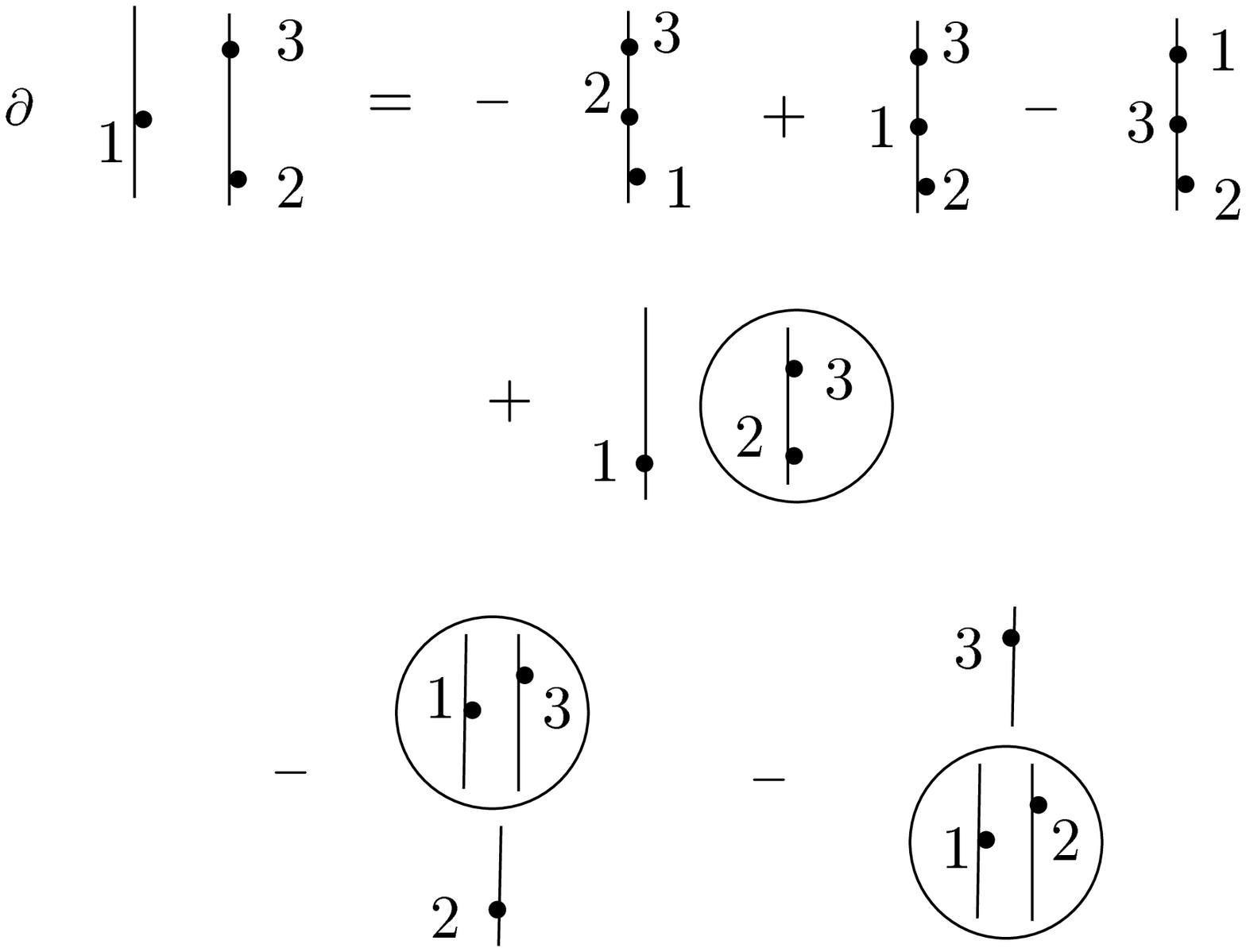}}}\fi
\begin{multline*}
\partial (\delta_3; 1| 2 3)
 \\
\begin{split}
=  & - (\delta_3; 123) + (\delta_3; 213) - (\delta_3; 231) \\
& + (\delta_2
\circ_2 \delta_2; (1|2) \circ_2 (12))\\
 & - \tau_{12} (\delta_2 \circ_2
\delta_2; (12) \circ_2 (1|2))
- (\delta_2 \circ_1
\delta_2; (12) \circ_1 (1| 2))
\end{split}
\end{multline*}
and similarly
\begin{multline*}
\partial (\delta_3; 12|3) \\
\begin{split}
 = & - (\delta_3; 123) + (\delta_3; 132) - (\delta_3; 312) \\
&- (\delta_2 \circ_1 \delta_2; (1|2)
\circ_1 (12))  \\
 & + (\delta_2 \circ_2 \delta_2; (12) \circ_2 (1|2))
+ \tau_{12} (\delta_2 \circ_1 \delta_2; (12) \circ_1 (1|2)),
\end{split}
\end{multline*}
where $\tau_{12}$ is the transposition exchanging 1
and 2, imply
\begin{multline}
\label{left}
d(\{v_1\}\{v_2,v_3\}) - \{dv_1\}\{v_2, v_3\}) \\
- (-1)^{\deg v_1 -1} \{
v_1\}\{dv_2, v_3\} - (-1)^{\deg v_1 + \deg v_2} \{v_1\}\{v_2, d v_3\}
\\
\begin{split}
 = & - M_3(v_1,v_2,v_3) + (-1)^{(\deg v_1) (\deg v_2)} M_3(v_2,v_1,v_3)
\\ & - (-1)^{\deg v_1 (\deg v_2 + \deg v_3)} M_3(v_2,v_3,v_1) + v_1 \circ ( v_2
\cdot v_3) \\
& - (-1)^{(\deg v_1 -1)\deg v_2} v_2 \cdot ( v_1 \circ
v_3) - (v_1 \circ v_2) \cdot v_3
\end{split}
\end{multline}
and
\begin{multline}
\label{right}
d(\{v_1,v_2\}\{v_3\}) - \{dv_1, v_2\} \{ v_3\}) \\
- (-1)^{\deg v_1 - 1} \{
v_1, dv_2 \} \{ v_3\} - (-1)^{\deg v_1 + \deg v_2} \{v_1, v_2\} \{ d
v_3\}
\\
\begin{split}
= & - M_3(v_1,v_2,v_3) + (-1)^{(\deg v_2) (\deg v_3)}
M_3(v_1,v_3,v_2) \\
& - (-1)^{\deg v_3 (\deg v_1 + \deg v_2)}
M_3(v_3,v_1,v_2)  - (v_1 \cdot v_2) \circ v_3 \\
& + v_1 \cdot ( v_2 \circ
v_3) + (-1)^{(\deg v_3 -1)\deg v_2} (v_1 \circ v_3) \cdot v_2 ,
\end{split}
\end{multline}
which can be regarded as the \emph{homotopy left and right Leibniz
rules} for the circle product with respect to the dot product, altered
by the trilinear $A_\infty$ product $M_3$.

Finally, for the 1-cell $(\delta_3; 123)$, we have

\let\picnaturalsize=N
\def\picsize{3.0in}
\def\picfilename{d32.eps}
\ifx\nopictures Y\else{\ifx\epsfloaded Y\else\input epsf \fi
\let\epsfloaded=Y
\centerline{\ifx\picnaturalsize N\epsfxsize \picsize\fi
\epsfbox{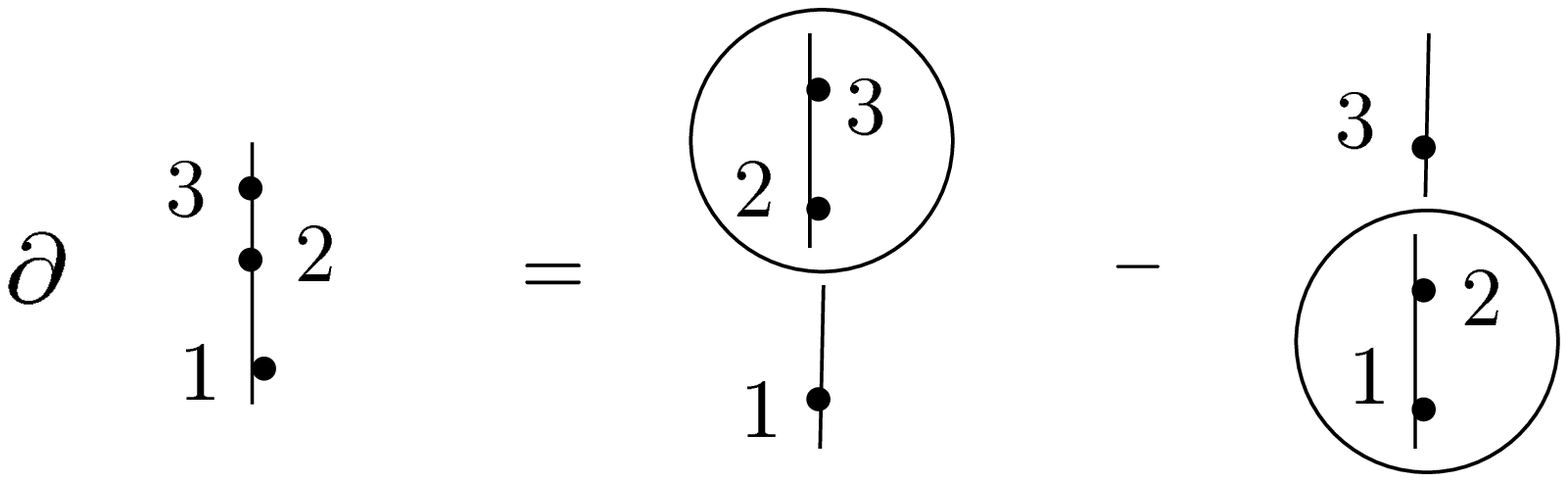}}}\fi
\[
\partial (\delta_3; 123) =
(\delta_2 \circ_1 \delta_2; (12)
\circ_1 (12)) - (\delta_2 \circ_2 \delta_2; (12) \circ_2 (12)),
\]
whence
\begin{multline}
\label{hass}
d(M_3 (v_1,v_2,v_3)) - M_3(dv_1, v_2, v_3) - (-1)^{\deg v_1} M_3 (v_1,
dv_2, v_3) \\
- (-1)^{\deg v_1 + \deg v_2} M_3 (v_1, v_2, d v_3)\\ =
(v_1 v_2) v_3 - v_1 (v_2 v_3),
\end{multline}
so that the dot product is \emph{homotopy associative} with $M_3$
being the \emph{associator} for it. In fact, the $M_n$'s form an
$A_\infty$-algebra, which is clearly seen at the operad level: the
$A_\infty$-operad, which is the natural cellular model of the
configuration operad of the real line, embeds into $K_\bullet \FF$
as a real part.

Note that the commutator $[\cdot,\cdot]$ of \eqref{bracket}
satisfies a graded \emph{homotopy Leibniz rule}:
\begin{multline}
\label{deriv}
[v_1 , v_2 v_3] - [v_1,v_2] v_3 - (-1)^{(\deg v_1 -1 ) \deg v_2} v_2
[v_1,v_3] \\
= - d ( \{v_1\}\{ v_2,v_3\} ) + \{dv_1\} \{ v_2, v_3\} +
(-1)^{\deg v_1 -1} \{v_1\} \{dv_2, v_3\} \\
+ (-1)^{\deg v_1 + \deg v_2} \{v_1\}
\{v_2, d v_3\} + (-1)^{\deg v_1 (\deg v_2 + \deg v_3) } (d ( \{
v_2,v_3\} \{v_1\} \\
- (-1)^{\deg v_2 + \deg v_3} \{ v_2, v_3\} \{dv_1\} -
\{dv_2, v_3\} \{v_1\} \\
- (-1)^{\deg v_2-1} \{v_2, d v_3\}\{v_1\}),
\end{multline}
which is a consequence of \eqref{left} and \eqref{right}.  The
\emph{homotopy Jacobi identity}
\begin{multline}
\label{jacobi}
[[v_1, v_2], v_3] + (-1)^{(\deg v_1-1) (\deg
v_2 + \deg v_3)} [[v_2, v_3], v_1]\\
+ (-1)^{(\deg v_3-1) (\deg v_1 +
\deg v_2)} [[v_3, v_1], v_2] \\
= \{\text{the differential of the sum of
all permutations of $\{v_1\}\{v_2\} \{v_3\}$}\}
\end{multline}
for the bracket $[,]$ follows from \eqref{distr}. More generally, the
graded symmetrizations of the operations $\{v_1\}\{v_2\} \dots
\{v_n\}$ form a homotopy Lie (or $L_\infty$-) algebra: the
symmetrizations are just the fundamental cycles of $\F{n}$, which
freely generate an operad with the differential coming from
contraction of edges of trees labeling a basis in this free operad;
see Beilinson-Ginzburg \cite{bg:1}. Due to Hinich-Schechtman
\cite{hs}, this is nothing but the homotopy Lie operad. In addition,
the circle product and the higher braces are a sort of universal
enveloping pre-Lie algebra, cf. \cite{gerst}, for the $L_\infty$-algebra.

\begin{rem}
The dot product $v_1 v_2$ and the bracket $[v_1,v_2]$ descend
to the cohomology of a $G_\infty$-algebra and endow the
cohomology with a G-algebra
structure; see Equations \eqref{hcomm}, \eqref{hass}, \eqref{deriv},
\eqref{jacobi}.
\end{rem}

\subsection{Homotopy Gerstenhaber, associative and Lie algebras and operads}
\label{hg-op}

While G-algebras combine commutative associative and Lie algebras, the
G-operad describing the class of G-algebras is a natural combination
of the commutative and the Lie operads $\comm$ and $\lie$, describing
the classes of commutative and Lie algebras, respectively. Indeed, we
can identify the latter operads with $H_0 (\FF,\nc)$ and $H_{n-1}
(\FF,\nc)$. The first identification follows from connectedness of the
spaces $\F{n}$: $H_0 (\F{n},\nc) \cong \nc \cong \comm (n)$, the
second follows Fred Cohen's Theorem \cite{fc}: $H_{n-1} (\F{n},\nc)
\cong H_{n-1} (F(n, \nc) , \nc) \cong \lie (n)$. The G-operad is just
$H_\bullet (\F{n}, \nc)$, according to another part of Cohen's Theorem, thereby
neatly interpolating the commutative and the Lie operads.

The associative operad $\ass$ can also be seen in this moduli space
picture: $H_0 (\FR{n},\nc) \cong H_0 (F(n,\nr), \nc) \cong \nc[S_n]
\cong \ass (n)$, where $F(n, \nr)$ is the configuration space of $n$
points on the real line, and $\FR{n}$ is the quotient of it by
translations and dilations, compactified similarly to $\F{n}$
above. The natural embedding $\nr \to \nc$, say, as the $y$ axis,
induces a morphism $\ass = H_0 (\FFR, \nc) \to H_\bullet (\FF, \nc)$
of the associative operad to the G-operad, of course through the
commutative one.

Similarly, as we have seen in the previous section,
$G_\infty$-algebras combine the structures of $A_\infty$- and
$L_\infty$-algebras. Accordingly, the $G_\infty$-operad neatly
incorporates the $A_\infty$- and the $L_\infty$-operads. The
$A_\infty$-operad is the operad consisting of all cells of the
real moduli space operad $\FFR$. It embeds as a cellular operad
into $K_\bullet \FF$ if we identify $\nr$ with the $y$ axis on
the complex plane $\nc$ and just consider the cells which are
formed by all points grouping on a vertical line. On the other
hand, the $L_\infty$-operad is the suboperad of the
$G_\infty$-operad $K_\bullet \FF$ generated by the sums of all
top cells of $\F{n}$'s, according to
Beilinson-Ginzburg-Hinich-Schechtman; see the end of the previous
section.

\section{Topological Field Theories}

The topological field theories described here are sometimes called {\sl
cohomological field theories} and are algebras over the operad of smooth
singular chains on moduli spaces of punctured Riemann spheres with
decorations.

\subsection{Reduced Conformal Field Theories}

Let $V$ be a topological vector space endowed with a collection of smooth
maps $\Psi:F(n,\nc) \,\to\, \Hom(V^{\otimes n},V)$, one for each $n\geq 1$,
taking $(z_1, z_2, \ldots, z_n)\linebreak[1] \mapsto \linebreak[0]
\Psi_{(z_1,z_2,\ldots, z_n)}$
satisfying the following axioms:
\begin{enumerate}
\item The map $\Psi$ should be $S_n$-equivariant, where the permutation group
$S_n$ acts on both $F(n,\nc)$ and $\Hom(V^{\otimes n},V)$ in the obvious ways.

\item $\Psi_{(a z_1+b,..., a z_n+b)} = \Psi_{(z_1,..., z_n)}$ for all $a \in
\nr^+$
and $b \in \nc$.

\item As the configuration of points $(z_1, ..., z_n)$ approaches a
composition of configurations $(t_1,..., t_p) \circ_i (w_1,..., w_{n-p+1})$
for some $i=1,\ldots, p$ then we have
\begin{multline*}
\Psi_{(z_1,..., z_n)}(v_1,..., v_n) \,
\longrightarrow\\ \Psi_{(t_1,..., t_p)}
(v_1,..., v_{i-1}, \Psi_{(w_1,..., w_{n-p+1})}(v_i,..., v_{i+n-p}),
v_{i+n-p+1},..., v_n),
\end{multline*}
where the compositions are understood to be in $\FF$. This algebraic
structure is said to be a {\sl reduced conformal field theory (reduced CFT).}
\end{enumerate}
The last axiom replaces the usual associativity axiom for VOAs and encodes
the various ways in which points can come together in the plane. The
operators $\Psi_{(z,0)} (v_1,v_2)$ are analogous to vertex
operators $Y(v_1,z)v_2$ of a VOA. It is clear
that we can redescribe this structure in the following way.
\begin{prop}
$V$ is a reduced conformal field theory if and only if $V$ is an
$\FF$-algebra.
\end{prop}

A reduced topological conformal field theory (reduced TCFT) can be
understood as an algebra over the smooth singular chain operad $C_\bullet
(\FF)$ of the configuration operad $\FF$ by analogy with TCFTs of
\cite{V}.  For technical reasons --- smooth singular chains, i.e., finite
sums of smooth mappings from simplices to the manifold, do not naturally
contain cells of $K_\bullet \FF$ among them --- we prefer to use the
following definition, imitating Segal's definition of a (full) TCFT
\cite{se}, see also the next section.

\begin{df}
\label{reduced}
A \emph{reduced topological conformal field theory $($reduced TCFT}) is
a complex $(V,d)$ of vector spaces, called a \emph{BRST complex or a
state space}, and a collection of operator-valued differential forms
$\Omega_n \in \Omega^\bullet (\F{n}, \Hom(V^{\otimes n},V))$, one for each
$n \ge 1$, such that
\begin{enumerate}
\item \[
\pi \Omega_n = \pi^* \Omega_n \qquad \text{for each $\pi \in S_n$},
\]
where $\pi$ on the left-hand side acts naturally on $\Hom(V^{\otimes n},V)$ and
$\pi^*$ is the geometric action of $\pi$ by relabeling the punctures,
\item \[
d_{\operatorname{DR}} \Omega_n = d_{\operatorname{BRST}} \Omega_n, \]
where $d_{\operatorname{DR}}$ is the de Rham differential
and $d_{\operatorname{BRST}}$ is the natural differential on
the  space $\Hom(V^{\otimes n},V)$ of $n$-linear operators on $V$,
\item
\[
\circ_i^* (\Omega_{m+n-1})  = \Omega_m \otimes \Omega_n \qquad
\text{ for each $i = 1, \dots, m$},
\]
where $\circ_i: \F{m} \times \F{n} \to \F{m+n -1}$ is the operad law.
\end{enumerate}
\end{df}

\begin{thm}
A reduced TCFT implies a natural $G_\infty$-algebra structure on the
BRST complex $V$.
\label{hga2}
\end{thm}
\begin{proof}
Given a reduced TCFT, we obtain the structure of an algebra over
the $G_\infty$-operad
\[
K_\bullet \F{n} \to \Hom (V^{\otimes n}, V), \qquad n \ge 1,
\]
by integrating the forms $\Omega_n$ over the cells:
\[
C \mapsto \int_C \Omega_n .
\]
\end{proof}

\subsection{Conformal Field Theories}
\label{tcft}

Let $\PP(n)$ be the moduli space of Riemann spheres with
$n+1$ distinct, ordered, holomorphically embedded unit disks which
do not overlap except, possibly, along their boundaries. The
permutation group on $n$ elements $S_n$ acts on $\PP(n)$ by permuting
the ordering of the first $n$ punctures. The collection $\PP =
\{\,\PP(n)\,\}$ forms an operad of complex manifolds where for all
$\Sigma$ in $\PP(n)$ and $\Sigma' $ in $\PP(n')$,
$\Sigma\circ_i\Sigma'$ in $\PP(n+n'-1)$ is obtained by cutting out the
$n'+1$st unit disk of $\Sigma'$, the $i$th disk of $\Sigma$ and
sewing along their boundaries using the identification $z\,\mapsto\,
\frac{1}{z}$. A {\sl (tree level, $c=0$) conformal field theory (CFT)}
is an algebra over $\PP$.  ``Tree level'' means that only genus 0
Riemann surfaces appear, while $c=0$ means that $V$ is a
representation of $\PP$ rather than a projective representation (see
Segal \cite{se} for details). We shall restrict to such CFTs for
simplicity. An important class of examples of holomorphic CFTs, that
is, CFTs where the algebra maps are holomorphic, come from vertex
operator algebras through the work of Huang-Lepowsky \cite{hl}.

\begin{df}[Segal \cite{se}, Getzler \cite{g}]
A \emph{topological conformal field theory $($TCFT}) is a complex $(V,d)$
and a collection of forms $\Omega_n \in \Omega^\bullet (\PP(n),
\Hom (V^{\otimes n}, V))$ satisfying the same axioms (1)--(3) of
Definition~\ref{reduced}. The complex $(V,d)$ is sometimes called
the \emph{BRST complex or the state space} of the theory in the
physics literature.
Also, the cohomology of the complex is called the \emph{space
of physical states}.
\end{df}
Again, we shall restrict to tree level, $c=0$
TCFTs, although physical examples of TCFT's do not always have $c=0$.
The forms $\Omega_n$ can be constructed using the operator formalism from
algebraic data known as a \emph{string background} (see \cite{agmv,ksv1}
for details). However, we shall not need this fact.

\begin{prop}
Let $(V,d)$ be a TCFT. Then $(V,d)$ is a reduced TCFT.
\end{prop}

\begin{proof}
This result follows from the fact that there exists a morphism of operads
(called {\sl string vertices}) $s:\FF\,\to \,\PP$ --- indeed, there exists
a
family of such morphisms. Consequently, any TCFT is a
reduced TCFT. A similar result was used by
Kimura, Stasheff, and Voronov in \cite{ksv1} to prove the existence of a
homotopy Lie algebra structure on a natural subcomplex of a TCFT.

The construction of $s$ uses a geometric result due to Wolf and Zwie\-bach
\cite{woz}, who showed that every conformal class of Riemann spheres with at
least three punctures can be endowed with a unique metric compatible with its
complex structure which solves a minimal area problem for metrics on the
punctured Riemann sphere subject to the constraint that the length of any
homotopically nontrivial closed curve on the punctured sphere is greater
than or equal to $2\pi$. This minimal area metric decomposes the punctured
sphere into flat cylinders foliated by closed geodesics with circumference
$2\pi$. A morphism $s$ is obtained for each pair $(l,f_l)$ where
$l$ is a real number greater than $\pi$ and $f_l:[2\pi,\infty)\,\to\,\nr$ is
a smooth monotonically increasing function such that $f_l(2\pi) = 2\pi$ and
$\lim_{x\,\to\,+\infty} f_l(x) = 2l.$

Let $\Sigma$ represent a point in $\F{n}$; since each irreducible
component of $\Sigma$ can be identified with the configuration of, say $m$,
points on $\nc$ quotiented by dilations and translations, assign to each
point a tangent direction which points along the positive real axis in $\nc$
and then assign the minimal area metric to each irreducible component of
$\Sigma$ minus its punctures and double points. Shrink all internal flat
cylinders in each irreducible component with a height $h$ greater than $2\pi$
to a flat cylinder with height $f_l(h)$. Furthermore, the minimal area metric
and the tangent directions at each puncture and on each side of every double
point gives rise to a holomorphically embedded unit disk centered there.
$s(\Sigma)$ in $\PP(n)$ is obtained by sewing together the irreducible
components on each side of every double point. The remaining curve has no
remaining double points and has a holomorphically embedded unit disk
around each puncture.
\end{proof}

Combining the previous proposition with Theorem~\ref{hga2}, we conclude
the following.

\begin{crl}
Let $(V,d)$ be a TCFT. Then $(V,d)$ admits the structure of a
$G_\infty$-algebra.
\end{crl}

The existence of an $A_\infty$-algebra structure on the state space of a TCFT
was observed in \cite{k}, cf. the homotopy associative structure of the
HIKKO open string-field theory group noticed by Stasheff \cite{jim:higher}.

\begin{rem}
In particular, we see (Section \ref{identities})
the structure of a homotopy commutative $A_\infty$-algebra
extending a dot product and an $L_\infty$-algebra extending a
bracket, naturally merged into one structure. This new
$L_\infty$ structure is not totally independent of the one
observed by Witten and Zwiebach \cite{wz,z} which was explained
operadically in \cite{ksv1}. The $L_\infty$ structure of the
present paper extends the one studied before. At the operadic
level, in the latter case, one projects along the phases
(that is why semirelative BRST cochains are needed), while in
the former one takes a section of this projection,
conveniently provided by the moduli space $\FF$ of punctures on
the sphere with an arrow at the $\infty$ puncture.
\end{rem}

\section{Applications of string vertices}

\subsection{Vassiliev knot invariants}
\label{vass}

Let $\Mh_{g,n}$ be the moduli space of stable curves of genus $g$ and $n$
ordered, distinct punctures whose double points are decorated with tangent
directions, one on each irreducible component on either side of each double
point, quotiented by the diagonal group of rotations by $U(1)$.
$\Mh_{g,n}$ is a $6g-6+2n$ dimensional compact, oriented orbifold with
corners which can be constructed by making real blowups along the
irreducible components of the divisor of the moduli space of stable curves
of genus $g$ and $n$ punctures $\Mc_{g,n}$. The collection of spaces $\Mh
= \{\,\Mh_{g,n}\,\}$ does not have natural composition maps between them,
unlike $\Mc = \{\,\Mc_{g,n}\,\}$, because for any two punctures which are
to be attached together, there is no natural way to choose tangent
directions at the double points. However, the space of (smooth) singular
chains $C_\bullet(\Mh) = \{\,C_\bullet(\Mh_{g,n})\,\}$ does have natural
composition maps between them by using the transfer which comes from
attaching two curves together at two punctures and then averaging over the
entire $S^1$ of tangent directions at the double points. $C_\bullet(\Mh)$
forms a generalization of an operad called a {\sl modular operad}, a
notion due to Getzler and Kapranov \cite{gek}, which generalizes operads
in two ways. The first is that there is no natural ``outgoing'' puncture
in this case --- any two punctures can be attached including two on a
single stable curve --- while the second is that higher genus stable
curves are allowed.  {\sl Throughout the remainder of this section, all
operads will be assumed to be modular operads unless otherwise stated.}

$\Mh_{g,n}$ is a stratified space whose strata are indexed by {\sl
stable $n$-graphs}. These are (connected) graphs with $n$ external
legs whose vertices are decorated with a nonnegative integer. One
associates a stable graph to each point in $\Mh_{g,n}$ by associating
to each irreducible component, a corolla with an external leg
associated to each puncture and double point on that component, an
integer assigned to that vertex corresponding to the genus of the
irreducible component, and then attaching the corollas together
whenever the irreducible components share a double point. A stratum
associated to a stable $n$-graph consists of all points in $\Mh_{g,n}$
whose associated graph is the given one. The stratification of
$\Mh_{g,n}$ is obtained by pulling up the stratification of
$\Mc_{g,n}$ via the canonical projection map $\Mh_{g,n} \,\to\,
\Mc_{g,n}$ which forgets the tangent directions at each double point.
The union of strata gives rise to a canonical filtration of $\Mh$
whose associated homology spectral sequence has an $E^1$-term $E^1 =
\{\,E^1_{g,n}\,\}$ which forms an operad of chain complexes. $E^1$
contains a suboperad $\G = \{\,\G_{g,n}\,\}$ called {\sl the graph
complex} of Kontsevich \cite{kon} which is nothing more than the
Feynman transform, a generalization of the operadic bar construction
to modular operads, of the commutative operad.

Kontsevich \cite{kon:vas} showed that the homology of the graph
complex $H_\bullet(\G_{g,n})$ has a special significance in the theory
of knots. He observed that the space of primitive chord diagrams of
order $n$ is isomorphic to $\oplus_{g=0}^n
H_0(\G_{g,n-g+1})_{S_{n-g+1}}$. The space of chord diagrams $\A$ then
is nothing more than the symmetric algebra over the space of primitive
chord diagrams, since $\A$ is a commutative, cocommutative Hopf
algebra. A {\sl weight system} is an element in $\A^*$. A key theorem
of Kontsevich \cite{kon:vas} and Bar-Natan \cite{bn} states that
weight systems are in one to one correspondence with finite type knot
invariants due to Vassiliev \cite{vas}. Such knot invariants as
the Jones, Alexander-Conway polynomials and their generalizations are
constructed from Vassiliev invariants. Most familiar examples of
weight systems are constructed from simple Lie algebras with an
invariant metric (see \cite{bn}) and their representations, although it
is has recently been discovered that not all weight systems are of
this kind \cite{vogel}. The following theorem was very likely known to
Kontsevich, who used homotopy Lie algebras with an invariant inner
product and graph complexes to construct Vassiliev knot invariants.

\begin{thm}
Let $(V,d)$ be an algebra over $\G$, Kontsevich's graph complex; then there
is an associated family of weight systems.
\end{thm}
\begin{proof}
The algebra structure morphism $\G\,\to\,\End{V}$ induces the morphism
$m:H_\bullet(\G)\,\to\,\End{H_\bullet(V)}$ which can be extended to a
homomorphism of (ungraded) commutative, associative algebras $m:\A\,\to\,\Sym
W$ where
$\Sym W$ denotes the (ungraded) symmetric algebra over the vector space $W$.
Here $W$ is a
$\nz$-graded vector space whose degree $n$ subspace is isomorphic to
$\oplus_{g=0}^n
\left(\End{H_\bullet(V),(g,n-g+1)}\right)_{S_{n-g+1}}$. Therefore, given a
vector in the symmetric algebra over $H_\bullet(V)$, one obtains a weight
system by contracting with a suitable inner product on $H_\bullet(V)$.
\end{proof}

Examples of such $\G$-algebras come from TCFT's modulo a slight technicality.
Let $(V,d)$ be a ($c=0$) TCFT. That is, let $\PP_{g,n}$ be the moduli space of
genus $g$ Riemann surfaces with $n$ holomorphically embedded unit disks which
do not intersect, except possibly along the boundaries in $\PP_{0,2}$. They
assemble into the operad $\PP = \{\,\PP_{g,n}\,\}$. A ($c=0$) TCFT is a
collection of endomorphism-valued differential forms $\Omega_{g,n}$, $n \ge 1$,
on $\PP_{g,n}$ satisfying the natural modular operad generalization of the
axioms of Definition~\ref{reduced} of Section~\ref{tcft}.

Let $U(1)$ be the subgroup of $\PP_{0,2}$ which consists of the set of
Riemann spheres with the standard chart about $0$ and the standard one around
$\infty$ rotated by multiplication by a phase. If $(V,d)$ is a TCFT then let
$\Delta:V\,\to\,V$ be the unary operation associated to $U(1)$ which is
regarded as a $1$-cocycle in $\PP_{0,2}$. The kernel $\Vr$ of $\Delta$
forms a subcomplex of $(V,d)$ called the {\sl semi-relative BRST
complex}.

Let $\Nc_{g,n}$ be the moduli space of stable curves of genus $g$ with $n$
distinct, ordered punctures which has the same decorations as points in
$\Mh_{g,n}$ but which have, in addition, tangent directions at each puncture.
The collection $\Nc = \{\,\Nc_{g,n}\,\}$ forms an operad. The string vertices
introduced in the previous section can also be described as a morphism of
(nonmodular) operads $s:\Nc\,\to\,\PP$ (string vertices) which are provided
by minimal area metrics. These minimal area metrics are proven to exist in
genus zero but are only conjectured to exist in higher genus \cite{woz}. We
shall assume that they exist in what follows.

Using the string vertices, we can pull back the forms $\Omega_{g,n}$ to
$\Nc_{g,n}$ and then push them forward to forms in $ \Omega^\bullet
(\Mh_{g,n}, \Hom(\Vr^{\otimes n}, \Vr))$, as in \cite{ksv1}.

\begin{thm}
Let $(V,d)$ be a TCFT; then $(\Vr,d)$ is an algebra over $\G$, thereby giving
rise to a family of weight systems.
\end{thm}
\begin{proof}
By integration of the forms $\Omega_{g,n}$, $(\Vr,d)$ becomes an algebra
over Kontsevich's graph complex $\G$ which appears in the top row in the
$E^1$ term in the homology spectral sequence of $\Mh$ associated to the
canonical filtration. Now apply the previous theorem.
\end{proof}

It is interesting that any two given TCFT's which are homotopic through the
space of TCFT's will give rise to isomorphic weight systems. Therefore,
there will be families of weight systems associated to each component of the
moduli space of TCFT's.

\subsection{Double loop spaces}
\label{double}

The following theorem, a byproduct of \emph{string vertices},
generalizes Stasheff's characterization of loop spaces as
$A_\infty$-spaces, i.e., algebras over his polyhedra operad, to double
loop spaces. It is also a refinement of the
Boardman-Vogt-May-Fadell-Neuwirth characterization of double loop
spaces considered up to homotopy as algebras over the little disks
operad $D(n)$, $n \ge 1$.

\begin{thm} Any double loop space is an algebra over the operad $\FF$. In
particular, the singular chain complex of a double loop space is a \hg
algebra, more precisely, an algebra over the singular chain operad
$C_\bullet (\FF)$. \end{thm}

\begin{proof}
The standard construction of Boardman and Vogt \cite{BoVo} provides a
double loop space with the natural structure of an algebra over the
little disks operad $D(n)$, $n \ge 1$. The same argument gives the
structure of an algebra over the operad $\PP$ of Riemann spheres with
holomorphic holes. String vertices deliver an operad morphism $\FF \to
\PP$, which yields a morphism $C_\bullet (\FF) \to C_\bullet (\PP)$.
\end{proof}

If we found a singular chain representative of each cell in $K_\bullet
(\FF)$ compatible with the operad structure, i.e., a morphism $K_\bullet
(\FF) \to C_\bullet (\FF)$ of operads, we would be able to answer the
following question. (In the case of Stasheff polyhedra, see more below, this
morphism does exist).
\begin{quest}
Is the singular chain complex of a double loop space a $G_\infty$-algebra?
\end{quest}
Note that string vertices offer an alternative approach to the study
of loop spaces compared to that given by May's approximation theory,
see for example J.~Stasheff's contribution \cite{jim:from} to this
volume. Ideally, approximation theory would provide a construction of
a space $K \FF X$ homotopy equivalent to the double loop space
$\Omega^2 \Sigma^2 X$ of the double suspension of a given topological
space $X$, such that $K\FF X$ is an algebra over the cellular operad
$K_\bullet \FF$. Another approach is the content of a promised theorem
of Getzler and Jones \cite[Introduction]{gj}.

It would be interesting to see whether string vertices exist in the
case of the little intervals operad. Here the collection of string
vertices should be nothing but a morphism of operads from the
compactified spaces of configurations of points on the real line to
the little intervals operad.  Since this configuration operad is
isomorphic to Stasheff's polyhedra operad, see Kontsevich \cite{kon},
it may yield a simpler ``quantum'' proof of Stasheff's famous theorem, saying
that any loop space is an algebra over the Stasheff polyhedra operad,
and in particular, the singular chain complex of a loop space is a
homotopy associative ($A_\infty$-) algebra.

\bibliographystyle{amsplain}

\begin{thebibliography}{10}

\bibitem{akman}
F.~Akman, \emph{On some generalizations of {B}atalin-{V}ilkovisky algebras},
  Preprint, Cornell University, 1995, {\tt{q-alg/9506027}}.

\bibitem{agmv}
L.~Alvarez-Gaume, C.~Gomez, G.~Moore, and C.~Vafa, \emph{Strings in the
  operator formalism}, Nuclear Phys. B \textbf{303} (1988), 455--521.

\bibitem{bn}
D.~Bar-Natan, \emph{On {V}assiliev knot invariants}, Topology \textbf{34}
  (1995), no.~2, 423--472.

\bibitem{bg:1}
A.~Beilinson and V.~Ginzburg, \emph{Infinitesimal structure of moduli spaces of
  {$G$}-bundles}, Internat. Math. Research Notices (1992), no.~4, 63--74.

\bibitem{BoVo}
J.~M. Boardman and R.~M. Vogt, \emph{Homotopy invariant algebraic structures on
  topological spaces}, Lecture Notes in Math., vol. 347, Springer-Verlag, 1973.

\bibitem{bor}
R.~E. Borcherds, \emph{Vertex operator algebras, {K}ac-{M}oody algebras and the
  {M}onster}, Proc. Natl. Acad. Sci. USA \textbf{83} (1986), 3068--3070.

\bibitem{C1}
F.~R. Cohen, \emph{The homology of {$\mathcal{C}_{n+1}$}-spaces, {$n\ge0$}},
  The homology of iterated loop spaces, Lecture Notes in Math., vol. 533,
  Springer-Verlag, 1976, pp.~207--351.

\bibitem{fc}
\bysame, \emph{Artin's braid groups, classical homotopy theory and sundry other
  curiosities}, Contemp. Math. \textbf{78} (1988), 167--206.

\bibitem{igor}
I.~B. Frenkel, Lectures at the {I}nstitute for {A}dvanced {S}tudy, January
  1988.

\bibitem{flm}
I.~B. Frenkel, J.~Lepowsky, and A.~Meurman, \emph{Vertex operator algebras and
  the {M}onster}, Academic Press, New York, 1988.

\bibitem{gerst}
M.~Gerstenhaber, \emph{The cohomology structure of an associative ring}, Ann.
  of Math. \textbf{78} (1963), 267--288.

\bibitem{gv1}
M.~Gerstenhaber and A.~A. Voronov, \emph{Higher order operations on
{H}ochschild complex}, Functional Anal. Appl. \textbf{29} (1995), no.~1, 1--6.

\bibitem{gv2}
\bysame, \emph{Homotopy {G}-algebras and moduli space operad}, Internat.\
  Math.\ Research Notices (1995), 141--153.

\bibitem{g}
E.~Getzler, \emph{{B}atalin-{V}ilkovisky algebras and two-dimensional
  topological field theories}, Commun. Math. Phys. \textbf{159} (1994),
  265--285, {\tt{hep-th/9212043}}.

\bibitem{gj}
E.~Getzler and J.~D.~S. Jones, \emph{Operads, homotopy algebra and iterated
  integrals for double loop spaces}, Preprint, Department of Mathematics, MIT;
  Department of Mathematics Northwestern University, March 1994,
  {\tt{hep-th/9403055}}.

\bibitem{gek}
E.~Getzler and M.~Kapranov, \emph{Modular operads}, Preprint, Department of
  Mathematics, MIT, August 1994, dg-ga/9408003.

\bibitem{gk}
V.~Ginzburg and M.~Kapranov, \emph{Koszul duality for operads}, Duke Math. J.
  \textbf{76} (1994), 203--272.

\bibitem{hs}
V.~Hinich and V.~Schechtman, \emph{Homotopy {L}ie algebras}, Adv. Studies Sov.
  Math. \textbf{16} (1993), 1--18.

\bibitem{h}
Y.-Z. Huang, \emph{Operadic formulation of topological vertex algebras and
  {G}erstenhaber or {B}atalin-{V}ilkovisky algebras}, Commun. Math. Phys.
  \textbf{164} (1994), 105--144, {\tt{hep-th/9306021}}.

\bibitem{hl}
Y.-Z. Huang and J.~Lepowsky, \emph{Vertex operator algebras and operads}, The
  Gelfand Mathematics Seminars, 1990--1992 (Boston), Birkh\"{a}user, 1993,
  {\tt{hep\-th/9301009}}, pp.~145--161.

\bibitem{k}
T.~Kimura, \emph{Operads of moduli spaces and algebraic structures in
  topological conformal field theory}, Moonshine, the Monster, and Related
  Topics (Providence) (C.~Dong and G.~Mason, eds.), Contemporary Math., vol.
  193, Amer. Math. Soc., 1996, pp.~159--190.

\bibitem{ksv1}
T.~Kimura, J.~Stasheff, and A.~A. Voronov, \emph{On operad structures of moduli
  spaces and string theory}, Commun. Math. Phys. \textbf{171} (1995), 1--25,
  {\tt{hep-th/9307114}}.

\bibitem{ksv2}
\bysame, \emph{Homology of moduli spaces of curves and commutative homotopy
  algebras}, The Gelfand Mathematics Seminars, 1993--1994 (J.~Lepowsky and
  M.~M. Smirnov, eds.), Birkh\"{a}user, 1996, to appear.

\bibitem{kon:sympl}
M.~Kontsevich, \emph{Formal (non)-commutative symplectic geometry}, The Gelfand
  Mathematics Seminars, 1990-1992 (L.~Corwin, I.~Gelfand, and J.~Lepowsky,
  eds.), Birkh\"{a}user, 1993, pp.~173--187.

\bibitem{kon:vas}
\bysame, \emph{Vassiliev's knot invariants}, Adv. Sov. Math. \textbf{16}
  (1993), no.~2, 137 -- 150.

\bibitem{kon}
\bysame, \emph{Feynman diagrams and low-dimensional topology}, First European
  Congress of Mathematics, Vol. II (Paris, 1992) (Basel), Progr. Math., vol.
  120, Birkh\"{a}user, 1994, pp.~97--121.

\bibitem{lz}
B.~H. Lian and G.~J. Zuckerman, \emph{New perspectives on the {BRST}-algebraic
  structure of string theory}, Commun. Math. Phys. \textbf{154} (1993),
  613--646, {\tt{hep-th/9211072}}.

\bibitem{loday}
J.-L. Loday, \emph{Une version non commutative des algebres de {L}ie: les
  algebres de {L}eibniz}, Enseign. Math. (2) \textbf{39} (1993), 269--293.

\bibitem{may:def}
J.~P. May, \emph{Definitions: Operads, algebras and modules}, Operads:
  Proceedings of Renaissance Conferences (J.-L. Loday, J.~Stasheff, and A.~A.
  Voronov, eds.), Amer. Math. Soc., 1996, in this volume, pp.~?--?

\bibitem{may}
\bysame, \emph{Operads, algebras and modules}, Operads: Proceedings of
  Renaissance Conferences (J.-L. Loday, J.~Stasheff, and A.~A. Voronov, eds.),
  Amer. Math. Soc., 1996, in this volume, pp.~?--?

\bibitem{nijenhuis}
A.~Nijenhuis, \emph{Jacobi-type identities for bilinear differential
  concomitants of certain tensor fields}, Indag. Math. \textbf{17} (1955),
  390--403.

\bibitem{se:old}
G.~Segal, \emph{Two-dimensional conformal field theories and modular functors},
  IXth Int. Congr. on Mathematical Physics (Bristol; Philadelphia) (B.~Simon,
  A.~Truman, and I.~M. Davies, eds.), IOP Publishing Ltd, 1989, pp.~22--37.

\bibitem{se:cam}
\bysame, Lectures at {C}ambridge {U}niversity, summer 1992.

\bibitem{se}
\bysame, \emph{Topology from the point of view of {Q.F.T}.}, Lectures at Yale
  University, March 1993.

\bibitem{jim}
J.~Stasheff, \emph{On the homotopy associativity of {H}-spaces, {II}}, Trans.
  Amer. Math. Soc. \textbf{108} (1963), 293--312.

\bibitem{jim:higher}
\bysame, \emph{Higher homotopy algebras: string field theory and {D}rinfeld's
  quasi-{H}opf algebras}, Proceedings of the XXth International Conference on
  Differential Geometric Methods in Theoretical Physics (New York, 1991) (River
  Edge, NJ), vol.~1, 1992, pp.~408--425.

\bibitem{jim:rec}
\bysame, \emph{Closed string field theory, strong homotopy {L}ie algebras and
  the operad actions of moduli space}, Perspectives on Mathematics and Physics
  (R.C. Penner and S.T. Yau, eds.), International Press, 1994,
  {\tt{hep-th/9304061}}, pp.~265--288.

\bibitem{jim:from}
\bysame, \emph{From operads to `phisically' inspired theories}, Operads:
  Proceedings of Renaissance Conferences (J.-L. Loday, J.~Stasheff, and A.~A.
  Voronov, eds.), Amer. Math. Soc., 1996, in this volume, pp.~?--?

\bibitem{vas}
V.~A. Vassiliev, \emph{Cohomology of knot spaces}, Theory of singularities and
  its applications (V.~I. Arnold, ed.), Amer. Math. Soc., 1990, pp.~23--69.

\bibitem{vogel}
P.~Vogel, \emph{Algebraic structures on modules of diagrams}, Preprint,
  Universit\'e Paris VII, August 1995.

\bibitem{V}
A.~A. Voronov, \emph{Topological field theories, string backgrounds and
  homotopy algebras}, Proceedings of the XXIInd International Conference on
  Differential Geometric Methods in Theoretical Physics, Ixtapa-Zihuatanejo,
  M\'exico (J.~Keller and Z.~Oziewicz, eds.), vol.~4, Advances in Applied
  Clifford Algebras (Proc. Suppl.), no.~S1, 1994, pp.~167--178.

\bibitem{wz}
E.~Witten and B.~Zwiebach, \emph{Algebraic structures and differential geometry
  in two-dimensional string theory}, Nucl. Phys. B \textbf{377} (1992),
  55--112.

\bibitem{woz}
M.~Wolf and B.~Zwiebach, \emph{The plumbing of minimal area surfaces}, Jour.
  Geom. Phys. \textbf{15} (1994), 23--56.

\bibitem{z}
B.~Zwiebach, \emph{Closed string field theory: {Q}uantum action and the
  {B}atalin-{V}ilkovisky master equa\-tion}, Nucl. Phys. B \textbf{390} (1993),
  33--152.

\end{thebibliography}

\providecommand{\bysame}{\leavevmode\hbox to3em{\hrulefill}\thinspace}

\end{document}